\documentclass[10pt,letterpaper,dvips]{article}
 \usepackage{amssymb,latexsym,amsmath,amsthm}
\newtheorem{thm*}{Theorem}

\newtheorem{thm}{Theorem}

\newtheorem{lemma}{Lemma}

\newtheorem{remark}{Remark}

\begin{document}

\def\d{ \partial } 
\def\Na{{\mathbb{N}}}

\def\Z{{\mathbb{Z}}}

\def\IR{{\mathbb{R}}}

\newcommand{\E}[0]{ \varepsilon}

\newcommand{\la}[0]{ \lambda}

\newcommand{\s}[0]{ \mathcal{S}}

\newcommand{\AO}[1]{\| #1 \| }

\newcommand{\BO}[2]{ \left( #1 , #2 \right) }

\newcommand{\CO}[2]{ \left\langle #1 , #2 \right\rangle} 

\newcommand{\R}[0]{ \IR\cup \{\infty \} } 

\newcommand{\co}[1]{ #1^{\prime}} 

\newcommand{\p}[0]{ p^{\prime}} 

\newcommand{\m}[1]{   \mathcal{ #1 }}

\newcommand{ \A}[1]{ \left\| #1 \right\|_H }

\newcommand{\B}[2]{ \left( #1 , #2 \right)_H }

\newcommand{\C}[2]{ \left\langle #1 , #2 \right\rangle_{  H^* , H } }

 \newcommand{\HON}[1]{ \| #1 \|_{ H^1} }

\newcommand{ \Om }{ \Omega}

\newcommand{ \pOm}{\partial \Omega}

\newcommand{\D}{ \mathcal{D} \left( \Omega \right)}

\newcommand{\DP}{ \mathcal{D}^{\prime} \left( \Omega \right)  }

\newcommand{\DPP}[2]{   \left\langle #1 , #2 \right\rangle_{  \mathcal{D}^{\prime}, \mathcal{D} }}

\newcommand{\PHH}[2]{    \left\langle #1 , #2 \right\rangle_{    \left(H^1 \right)^*  ,  H^1   }    }

\newcommand{\PHO}[2]{  \left\langle #1 , #2 \right\rangle_{  H^{-1}  , H_0^1  }} 

 \newcommand{\HO}{ H^1 \left( \Omega \right)}

\newcommand{\HOO}{ H_0^1 \left( \Omega \right) }

\newcommand{\CC}{C_c^\infty\left(\Omega \right) }

\newcommand{\N}[1]{ \left\| #1\right\|_{ H_0^1  }  }

\newcommand{\IN}[2]{ \left(#1,#2\right)_{  H_0^1} }

\newcommand{\INI}[2]{ \left( #1 ,#2 \right)_ { H^1}} 

\newcommand{\HH}{   H^1 \left( \Omega \right)^* } 

\newcommand{\HL}{ H^{-1} \left( \Omega \right) }

\newcommand{\HS}[1]{ \| #1 \|_{H^*}}

\newcommand{\HSI}[2]{ \left( #1 , #2 \right)_{ H^*}}

\newcommand{\WO}{ W_0^{1,p}} 
\newcommand{\w}[1]{ \| #1 \|_{W_0^{1,p}}}  

\newcommand{\ww}{(W_0^{1,p})^*}

\title{Regularity of the extremal solution  in a MEMS model with advection}

\author{\sc{Craig COWAN\footnote{Department of Mathematics, University of British Columbia, Vancouver, B.C. Canada V6T 1Z2. E-mail: cowan@math.ubc.ca. }
\quad
\quad
Nassif GHOUSSOUB\footnote{Department of Mathematics, University of British Columbia, Vancouver, B.C. Canada V6T 1Z2. E-mail: nassif@math.ubc.ca. Research partially supported by the Natural Science and Engineering Research Council of Canada.}\quad
}}
\date{\today}

\smallbreak
\maketitle

\begin{abstract} We consider the regularity of the extremal solution of the nonlinear eigenvalue problem 
 \begin{equation*} 
(S)_\lambda \qquad \left\{ 
\begin{array}{rcr}
-\Delta u + c(x) \cdot \nabla u &=& \frac{\lambda}{(1-u)^2} \qquad \mbox{ in $ \Omega$}, \\
u &=& 0 \qquad \mbox{ on $ \pOm$}, \\
\end{array}
\right.
\end{equation*} where $ \Omega $ is a smooth bounded domain in $ \IR^N$ and $ c(x)$ is a smooth bounded vector field on $\bar \Omega$.   We show that, just like in the advection-free model ($c\equiv 0$), 
all semi-stable solutions are smooth if (and only if) the dimension $N\leq 7$. The novelty here comes from the lack of a suitable variational characterization for the semi-stability assumption.  We overcome this difficulty by using  a general  version of Hardy's inequality.  In a forthcoming paper \cite{CG2}, we indicate how this method applies to many other nonlinear eigenvalue problems involving advection (including the Gelfand problem), showing that they all essentially have the same critical dimension as their advection-free counterparts. 

 \end{abstract}

\section{Introduction} 

The following equation has often been used to model a simple \emph{Micro-Electro-Mechanical System} (MEMS) device:
\begin{equation*} 
(P)_\lambda \qquad \left\{ 
\begin{array}{rcr}
-\Delta u &=& \frac{\lambda}{(1-u)^2} \qquad \mbox{ in $ \Omega$}, \\
u &=& 0 \qquad  \mbox{\quad  on $ \pOm$}, \\
\end{array}
\right.
\end{equation*} where $ \Omega $ is a smooth bounded domain in $ \IR^N$, $ \lambda >0$ is proportional to the applied voltage and $ 0 < u(x) <1$ denotes the deflection of the membrane.   This model has been extensively studied, see \cite{P}, \cite{BP} in regards to the model and \cite{GG}, \cite{EGG}, \cite{GPW} for mathematical aspects of $(P)_\lambda$.     It is well known (see above references) that there exists some positive finite critical parameter $ \lambda^*$ such that for all $ 0 < \lambda < \lambda^*$,  the equation $(P)_\lambda$ has a smooth minimal stable (see below) solution  $ u_\lambda$, while for $ \lambda > \lambda^*$ there are no weak solutions of $(P)_\lambda$ (see \cite{GG} for a precise definition of weak solution).    Standard elliptic regularity theory yields that a solution $u$ of $ (P)_\lambda$ is smooth if and only if $ \sup_\Omega u <1$.    One can also show that $ \lambda \mapsto u_\lambda(x)$ is increasing and hence one can define the extremal solution 
\[ u^*(x):= \lim_{\lambda \nearrow \lambda^*} u_\lambda(x), \] which can be shown to be a weak solution of $(P)_{\lambda^*}$.     

Recall that a smooth solution $u$ of $ (P)_\lambda$ is said to be \emph{minimal} if  any other solution $ v$ of $ (P)_\lambda$ satisfies $ u \le v $ a.e. in $ \Omega$.   Such solutions are then  \emph{semi-stable} meaning that the principal eigenvalue of the linearized operator 
\[ L_{u,\lambda}:= -\Delta - \frac{2 \lambda}{(1-u)^3} \]  in $H_0^1(\Omega)$ is nonnegative.  This property can be expressed variationally by  the inequality
\begin{equation} \label{stability-var}
2 \lambda \int_\Omega \frac{ \psi^2}{(1-u)^3} \le \int_\Omega | \nabla \psi|^2, \qquad \forall \psi \in H_0^1(\Omega).
\end{equation} which can be viewed as the nonnegativeness  of the second variation of the  energy functional associated with $ (P)_\lambda$ at $ u$. 


Now a question of interest is whether $ u^*$ is a smooth solution of $ (P)_{\lambda^*}$.    It is shown in \cite{GG} that this is indeed the case provided $ N \le 7$.  This result is optimal in the sense that $ u^*$ is singular in dimension $ N \ge 8$ with $ \Omega$ taken to be the unit ball.   \\

Our main interest here will be in the regularity of the extremal solution associated with 
\begin{equation*} 
(S)_\lambda \qquad \left\{ 
\begin{array}{rcr}
-\Delta u + c(x) \cdot \nabla u &=& \frac{\lambda}{(1-u)^2} \qquad \mbox{ in $ \Omega$}, \\
u &=& 0 \qquad \mbox{ on $ \pOm$}, \\
\end{array}
\right.
\end{equation*} where $ c \in C^\infty( \overline{ \Omega}, \IR^N)$ and where again $ \Omega$  is a smooth bounded domain in $ \IR^N$.  Modifying the proofs used in analyzing $(P)_\lambda$ one can again show  the existence of a positive finite critical parameter $ \lambda^*$ such that for $ 0< \lambda < \lambda^*$ there exists a smooth minimal solution $u_\lambda$ of $(S)_\lambda$, while  there are no smooth solutions of $(S)_\lambda$ for $ \lambda > \lambda^*$.  Moreover, the minimal solutions are also {\it semi-stable} in the sense that the principal eigenvalue of the corresponding linearized operator 
\[ 
L_{u,\lambda, c}:= -\Delta + c(x) \cdot \nabla  - \frac{2 \lambda}{(1-u_\lambda)^3}
 \]  in $H_0^1(\Omega)$ is non-negative.
See \cite{explosion} where these results are proved for general $C^1$ convex nonlinearities which are  superlinear at $ \infty$.     
Our main result concerns the regularity of the extremal solution of $(S)_\lambda$. 

\begin{thm} \label{main} If $1\leq N \le 7$, then the extremal solution $ u^*$ of $(S)_{\lambda^*}$ is smooth. 

\end{thm}   

\begin{remark} \rm A crucial (in fact the main) ingredient  in proving the regularity of $u^*$ in $(P)_\lambda$,  is the energy inequality (\ref{stability-var}) which is used in conjunction with the equation $(P)_\lambda$, to obtain uniform (in $ \lambda$) $L^p$-estimates on $ (1-u_\lambda)^{-2}$ whenever $u_\lambda$ is the minimal solution (See \cite{GG}).   However,   the semi-stability of $ u_\lambda$ in the case of  $(S)_\lambda$,
does not  translate into an \emph{energy inequality} which allows the use of arbitrary test functions.  Overcoming this will be the major hurdle in proving Theorem \ref{main}.  

We point out, however, that if $ c(x) = \nabla \gamma$ for some smooth function $\gamma$ on $\bar \Omega$,  then the semi-stability condition on the minimal solution $u_\lambda$ of $ (S)_\lambda$ translates into
\begin{equation} \label{stability-var.2}
2 \lambda \int_\Omega \frac{ e^{-\gamma}\psi^2}{(1-u_\lambda)^3} \le \int_\Omega e^{-\gamma}| \nabla \psi|^2, \qquad \forall \psi \in H_0^1(\Omega).
\end{equation} 
Then, with slight modifications, one can use the standard approach for $(P)_\lambda$ to obtain the analogous result for $(S)_\lambda$ stated in Theorem \ref{main}.

The novel case is therefore when $c$ is a divergence free vector field. Actually, we shall use the following version of the Hodge decomposition, in order to deal with general vector fields $c$.  


\begin{lemma} Any vector field  $ c \in C^\infty(\overline{\Omega},\IR^N)$ can be decomposed as $ c(x) = - \nabla \gamma + a(x)$ where $ \gamma $ is a smooth scalar function and $ a(x)$ is a smooth bounded vector field such that  $ {\rm div}( e^\gamma a)=0$.   

\end{lemma} 

\begin{proof}   By the Krein-Rutman theory,    the linear eigenvalue problem 
\begin{equation} \label{lin.0}
\left\{
\begin{array}{rcl}
\Delta \alpha + {\rm div}( \alpha c) &=& \mu \alpha \qquad \Omega,  \\
(\nabla \alpha +\alpha c)\cdot n&=& 0 \qquad  \quad  \pOm,
\end{array}
\right.
\end{equation}
where $n$ is the unit outer normal on $\partial \Omega$, has a positive solution $\alpha$ in $\Omega$ when $\mu$  is the principal eigenvalue. Integrating the equation over $\Omega$, one sees that $\mu=0$. The positivity of $\alpha$ on the boundary follows from the boundary condition and the maximum principle. In other words, we have that $\Delta \alpha + {\rm div}( \alpha c) =0$ on $\Omega$, and $\alpha>0$ on $\bar \Omega$. 

  Now define $ \gamma:=\log( \alpha)$ and $ a:= c + \nabla \gamma$.   An easy computation shows that $ {\rm div}( e^\gamma a)=0$. 
\end{proof}  
Throughout the rest of this note $ c, a, \gamma $ will be defined as above.

\end{remark} 

\section{A general Hardy inequality and non-selfadjoint eigenvalue problems} 

Consider the linear eigenvalue problem
\begin{equation} \label{lin}
\left\{
\begin{array}{rcl}
-\Delta \phi + c \cdot \nabla \phi -\rho \phi&=& K \phi \qquad \Omega,  \\
\phi &=& 0 \qquad  \quad  \pOm,
\end{array}
\right.
\end{equation}  
where $c$ is a smooth bounded vector field on $\Omega$, $\rho \in C^\infty(\overline{\Omega})$ and $K$ is a scalar. We assume that 
 $ (\phi,K)$ is the principal eigenpair for (\ref{lin}) and that  $ \phi>0$ in $ \Omega$, and $K \geq 0$.  Note that elliptic regularity theory shows that $ \phi$ is then smooth.

We shall now use a general Hardy inequality to make up for the lack of a variational characterization for the pair $ (\phi,K)$.   The following result is taken from  \cite{craig}, which we duplicate here for the convenience of the reader.  For a complete discussion on general Hardy inequalities including best constants, attainability and improvements of, see \cite{craig}.   We should point out that this approach to Hardy inequalities is not new, but it is generally restricted to specific functions $E$ which yield known versions of Hardy inequalities; see \cite{AS} and reference within.

\begin{lemma}  Let $ A(x)$ denote a uniformly positive definite $ N \times N$ matrix with smooth coefficients defined on $\Omega$.  Suppose $ E$ is a smooth positive function on $ \Omega$ and fix a constant $\beta$ with $ 1 \le \beta \le 2$.    Then,   for all $\psi \in H_0^1(\Omega)$ we have 
\begin{equation} \label{ppp}
\int_\Omega | \nabla \psi|_A^2 \ge \frac{ \beta (2-\beta)}{4} \int_\Omega \frac{ | \nabla E|_A^2}{E^2} \psi^2 + \frac{\beta}{2} \int_\Omega \frac{ -{\rm div}(A \nabla E)}{E} \psi^2,
\end{equation} where $ \int_\Omega | \nabla \psi|_A^2 = \int_\Omega A(x) \nabla \psi \cdot \nabla \psi$.
\end{lemma} 

\begin{proof} For simplicity we prove the case where  $ A(x) $ is given by the identity matrix. For the general case, we refer to \cite{craig}.  Let $ E_0$ denote a smooth positive function defined in $ \Omega$ and let $ \psi \in C_c^\infty(\Omega)$.   Set $ v:= \frac{ \psi}{\sqrt{E_0}}$.   Then 
\begin{equation} \label{zz} | \nabla \psi|^2 = E_0 | \nabla v|^2 + \frac{ | \nabla E_0|^2 }{4 E_0^2} \psi^2 + v \nabla v \cdot \nabla E_0.
\end{equation}
  Integrating the last term by parts gives 
\[
  \int_\Omega v \nabla v \cdot \nabla E_0 = \frac{1}{2} \int_\Omega \frac{ -\Delta E_0}{E_0} \psi^2 \] and so integrating (\ref{zz}) gives 
\begin{equation} \label{zzz}
 \int_\Omega | \nabla \psi|^2 \ge \frac{1}{4} \int_\Omega \frac{ | \nabla E_0|^2}{E_0^2} \psi^2 + \frac{1}{2} \int_\Omega \frac{-\Delta E_0}{E_0} \psi^2,
 \end{equation}  where we dropped a nonnegative term.    So we have the desired result for $ \beta=1$.   When $ \beta \neq 1$ one puts $ E_0:=E^\beta$ into (\ref{zzz}) and collects like terms to obtain the desired result. 
\end{proof}    

We
now use the above lemma to obtain an energy inequality valid for the principal eigenpair of (\ref{lin}). 

\begin{thm} \label{nsa.var} Suppose that the principal eigenpair $(\phi,K)$ of (\ref{lin}) are such that $ \phi>0$ and $K\geq 0$.   Then,  for $ 1 \le \beta \le 2$ we have for all $ \psi \in H_0^1(\Omega)$,  
\begin{equation} \label{po}
\int_\Omega e^\gamma | \nabla \psi|^2 \ge \frac{\beta(2-\beta)}{4} \int_\Omega \frac{ e^\gamma | \nabla \phi|^2}{\phi^2} \psi^2 + \frac{\beta}{2} \int_\Omega   e^\gamma \rho (x) \psi^2 - \frac{\beta}{2} \int_\Omega \frac{e^\gamma a \cdot \nabla \phi}{\phi} \psi^2. 
\end{equation} 
\end{thm}

\begin{proof} Note that (\ref{lin}) can be rewritten as 
\[ - {\rm div}( e^\gamma \nabla \phi) + e^\gamma a \cdot \nabla \phi = e^\gamma \left(  \rho(x) +K \right) \phi \qquad \mbox{ in $ \Omega$}, \]  where as mentioned above we are using the decomposition $ c = - \nabla \gamma + a$.    We now set $ E:= \phi$ and  $ A(x) = e^\gamma I$ (where $I$ is the identity matrix) and use (\ref{ppp}) along with the above equation to obtain the desired result. Note that we have dropped the nonnegative term involving $K$.   

\end{proof}

 \section{Proof of theorem \ref{main}}  
 
For $ 0 < \lambda < \lambda^*$, we   denote by $ u_\lambda $  the smooth minimal semi-stable solution of  $(S)_\lambda$.  Let $ (\phi,K)$ denote the principal eigenpair associated with the linearization of $(S)_\lambda$ at $ u_\lambda$.   Then $ 0 < \phi$ in $ \Omega$, $ 0 \le K$ and $ (\phi,K)$ satisfy
\begin{equation} \label{non.lin}
\left\{
\begin{array}{rcl}
-\Delta \phi + c \cdot \nabla \phi &=& ( \frac{2 \lambda}{(1-u_\lambda)^3}+K) \phi \qquad \Omega,  \\
\phi &=& 0 \qquad \qquad \qquad \qquad \pOm.
\end{array}
\right.
\end{equation}
Again, elliptic regularity theory shows that $ \phi$ is smooth.  Consider  $ c=-\nabla \gamma + a$ to be the decomposition of $c$ described in Lemma 1. We now obtain the main estimate.  

\begin{thm} For $ 0 < \lambda < \lambda^*$, $ 1 < \beta <2$ and $ 0<t< \beta + \sqrt{ \beta^2+\beta}$, we have the following estimate: 
\[ \lambda \left( \beta - \frac{t^2}{2t+1} \right) \int_\Omega \frac{e^\gamma}{(1-u_\lambda)^{2t+3}} \le 2 \beta \lambda \int_\Omega \frac{e^\gamma}{(1-u_\lambda)^{t+3}} + \frac{ \beta \| a \|_{L^\infty}^2}{4 (2-\beta)} \int_\Omega \frac{ e^\gamma}{(1-u_\lambda)^{2t}}. \] 
\end{thm}  

\begin{proof}  Fix $ 0 < \beta <2$, let $ 0 <t$ and $ u $ denote the minimal solution associated with $(S)_\lambda$. We shall use Theorem \ref{nsa.var} with $\rho (x)=\frac{2 \lambda}{(1-u_\lambda)^3}$.  
 Put $ \psi:= \frac{1}{(1-u)^t}-1$ into (\ref{po}) to obtain 
\begin{eqnarray*}
t^2 \int_\Omega \frac{ e^\gamma | \nabla u|^2}{(1-u)^{2t+2}} & \ge & \beta \lambda \int_\Omega \frac{e^\gamma}{(1-u)^3} \left( \frac{1}{(1-u)^t}-1 \right)^2  \\
&&+ \frac{\beta}{2} \int_\Omega e^\gamma \left( \frac{(2-\beta)}{2} \frac{ | \nabla \phi|^2}{\phi^2} - \frac{a \cdot \nabla \phi}{\phi} \right) \psi^2.   
\end{eqnarray*} Now note that $(S)_\lambda $ can be rewritten as 
\[ -{\rm div}( e^\gamma \nabla u) + e^\gamma a \cdot \nabla u = \frac{ \lambda e^\gamma}{(1-u)^2} \qquad \mbox{ in $ \Omega$},\] and test this on $ \bar{ \phi}:= \frac{1}{(1-u)^{2t+1}}-1 $ to obtain  
\[ (2t+1) \int_\Omega \frac{ e^\gamma | \nabla u|^2}{(1-u)^{2t+2}} + H = \lambda \int_\Omega \frac{ e^\gamma}{(1-u)^2} \left( \frac{1}{(1-u)^{2t+1}}-1 \right),  \]  where 
\[ H:=\int_\Omega e^\gamma a \cdot \nabla u \left( \frac{1}{(1-u)^{2t+1}}-1 \right).\] One easily sees that $H=0$ after considering the fact $H$ can be rewritten in the form $ \int_\Omega ( e^\gamma a) \cdot \nabla G(u) $ for an appropriately chosen function $ G $ with $ G(0)=0$. 
  Combining the above two inequalities and dropping some positive terms gives 
\begin{eqnarray*}
\lambda \left( \beta - \frac{t^2}{2t+1} \right) \int_\Omega \frac{e^\gamma}{(1-u)^{2t+3}} & \le & 2 \beta \lambda \int_\Omega \frac{e^\gamma}{(1-u)^{3+t}} \\
&& + \frac{\beta}{2} \int_\Omega e^\gamma \Lambda(x) \left( \frac{1}{(1-u)^t} -1 \right)^2
\end{eqnarray*} where 
\[ \Lambda(x):= \frac{a \cdot \nabla \phi}{\phi} - \frac{(2 - \beta)}{2} \frac{ |\nabla \phi|^2}{\phi^2}.\]   
Simple calculus shows that 
\[ \sup_\Omega \Lambda(x) \le \frac{ \| a \|_{L^\infty}^2}{2(2-\beta)},\]  which, after substituting into the above inequality, completes the proof of the main estimate. \end{proof}

Note now that  the restriction $ t< \beta + \sqrt{ \beta^2+\beta} $ is needed to ensure that the coefficient $ \beta - \frac{t^2}{2t+1} $ is positive. It follows then that   
  $ \frac{1}{(1-u_\lambda)^2} $ is uniformly  bounded (in $ \lambda$) in $ L^p(\Omega)$ for all $ p< p_0:= \frac{7}{2} + \sqrt{6} \approx 5.94...$ and after passing to limits we have the same result for the extremal solution $ u^*$.  
  
  To conclude the proof of Theorem \ref{main}, it suffices to note the following result. 

\begin{lemma} Suppose $ 3 \le N \le 7$ and the extremal solution $ u^* $ satisfies $ \frac{1}{(1-u^*)^2} \in L^\frac{3N}{4} (\Omega)$, then $ u^*$ is smooth. 

\end{lemma} 

\begin{proof}  First note that by elliptic regularity one has $ u^* \in W^{2, \frac{3N}{4}} (\Omega)$ and after applying the Sobolev embedding theorem one has $ u^* \in C^{0,\frac{2}{3}} ( \overline{ \Omega})$.  Now suppose $ \| u \|_{L^\infty}=1$ so that there is some $ x_0 \in \Omega $ such that $ u(x_0)=1$.   Then 
\[ \frac{1}{1-u(x)} \ge \frac{C}{|x-x_0|^\frac{2}{3}} \] and hence  
\[ \infty > \int_\Omega \frac{1}{ ((1-u^*)^2)^\frac{3N}{4}} \ge C \int_\Omega \frac{1}{|x|^N} = \infty, 
\] 
which is a contradiction. It follows that $ \frac{1}{(1-u^*)^2} \in L^\infty (\Omega)$, and   $ u^*$ is therefore smooth.
\end{proof}   
Using this lemma and the above $L^p$-bound on $ \frac{1}{(1-u^*)^2}$, one sees that $ u^*$ is smooth for $ 3 \le N \le 7$.   To show the result in dimensions $N=1,2$, one needs a slight variation of the above argument. We omit the details, and the interested reader can consult \cite{GG} for the proof when $ c(x)=0$. 

\begin{remark} \rm As mentioned in the abstract, this method applies to most non-selfadjoint eigenvalue problems of the form  
\begin{equation*} 
(S)_\lambda \qquad \left\{ 
\begin{array}{rcr}
-\Delta u + c(x) \cdot \nabla u &=& \lambda f(u) \qquad \mbox{ in $ \Omega$}, \\
u &=& 0 \qquad \mbox{ on $ \pOm$}, \\
\end{array}
\right.
\end{equation*}
where $f(u)$ is an appropriate convex nonlinearity such as $f(u)=e^u$, and $f(u)=(1+u)^p$. It shows in particular that the presence of an advection does not change the critical dimension of the problem, hence addressing an issue raised recently by Berestycki et al \cite{explosion}.     One can also extend the general regularity results of Nedev \cite{Ned} (for general convex $f$ in dimensions $2$ and $3$) and those of Cabre and Capella \cite{CaCa1}  (for general radially symmetric $f$ on a ball, and up to dimension $9$). All these questions are the subject of a forthcoming paper \cite{CG2}.
\end{remark}

\end{document}